\documentclass{birkjour}
\usepackage{amsmath}

\newtheorem{thm}{Theorem}[section]

\theoremstyle{definition}

\theoremstyle{remark}

\numberwithin{equation}{section}

\begin{document}

\title[Weingarten and Linear Weingarten Canal Surfaces]{Weingarten and
Linear Weingarten \\
Canal Surfaces}
\author[Y\i lmaz TUN\c{C}ER and Dae Won YOON]{Y\i lmaz TUN\c{C}ER and Dae
Won YOON}
\email{yilmaz.tuncer@usak.edu.tr}
\thanks{}
\email{dwyoon@gnu.ac.kr}
\date{January 6, 2015}
\dedicatory{}

\begin{abstract}
In this study, we investigated the (K,H),(K,K$_{\text{II}}$), (H,K$_{\text{II%
}}$)-Weingarten and (K,H),(K,K$_{\text{II}}$),(H,K$_{\text{II}}$) and (K,H,K$%
_{\text{II}}$)-linear Weingarten canal surfaces in IR$^{3}$.
\end{abstract}

\maketitle

%
%
%
%
%
%
%
%
%


\address{U\c{s}ak University \\ Science and Art Faculty \\
Mathematics Department \\ U\c{s}ak TURKEY}


\address{Department of Mathematic Education and
RINS \\ Gyeongsang National University \\ Jinju 660- 701 \\ Repuclic of Korea}





\section{Introduction}

In 1863, Julius Weingarten was able to make a major step forward in the
topic when he gave a class of surfaces isometric to a given surface of
revolution. Surface for which there is a definite functional relation
between the principal curvatures (which called curvature diagram) and also
between the Gaussian and the mean curvatures is called Weingarten surface.
The knowledge of first fundamental form I and second fundamental form II of
a surface facilitates the analysis and the classification of surface shape.
Especially recent years, the geometry of the second fundamental form II has
become an important issue in terms of investigating intrinsic and extrinsic
geometric properties of the surfaces. Very recent results concerning the
curvature properties associated to II and other variational aspects can be
found in \cite{Sv}. One may associate to such a surface M geometrical
objects measured by means of its second fundamental form, as second Gaussian
curvature K$_{\text{II}}$, respectively. We are able to compute K$_{\text{II}%
}$ of a surface by replacing the components of the first fundamental form $%
E, $ $F,$ $G$ by the components of the second fundamental form $e,$ $f,$ $g$
respectively in Brioschi formula which is given by Francesco Brioschi in the
years of 1800's. Identification of the curvatures related to the second
fundamental form of a surface opened a door to research the new classes of
Weingarten surfaces. Since the middle of the last century, several geometers
have studied Weingarten surfaces and linear Weingarten surfaces and obtained
many interesting and valuable results\cite{3,4,7,12,14,16}. For study of
these surfaces, W. K\"{u}hnel and G.Stamou investigate ruled $($X,Y$)$%
-Weingarten surface in Euclidean 3-space E$^{3}$\cite{12,16}.
Also,C.Baikoussis and Th. Koufogiorgos studied helicoidal (H,K$_{\text{II}}$%
)-Weingarten surfaces\cite{1}. F.Dillen and W. K\"{u}hnel and F.Dillen and
W.Sodsiri gave a classification of ruled $($X,Y$)$-Weingarten surface in
Minkowski 3-space E$_{1}^{3}$, where X,Y$\in \left\{ \text{K,H,K}%
_{II}\right\} $\cite{3,4,5}. D. Koutroufiotis and Th.Koufogiorgos and T.
Hasanis investigate closed ovaloid $($X,Y$)$-linear Weingarten surface in E$%
^{3}$\cite{10,11}. D. W. Yoon and D.E.Blair and Th.Koufogiorgos classified
ruled $($X,Y$)$-linear Weingarten surface in E$^{3}$\cite{2,17}. D. W. Yoon
and J.S.Ro studied tubes in Euclidean 3-space which are (K,H),(K,K$_{\text{II%
}}$),(H,K$_{\text{II}}$)-Weingarten and linear Weingarten tubes\cite{15}. D.
W. Yoon also studied the Weingarten and linear Weingarten types translation
surfaces in Euclidean 3-space.

Surface theory has been a popular topic for many researchers in many
aspects. Besides the using curves and surfaces, canal surfaces are the most
popular in computer aided geometric design such that designing models of
internal and external organs, preparing of terrain-infrastructures,
constructing of blending surfaces, reconstructing of shape, robotic path
planning, etc. (see, \cite{Fr,Wm,Xu}).

In this study, we investigated the $($K,H$),$ $($K,K$_{\text{II}}),$ $($H,K$%
_{\text{II}})-$Weingarten and $($K,H$),($K,K$_{\text{II}}),($H,K$_{\text{II}%
})$ and $\left( \text{K,H,K}_{\text{II}}\right) -$linear Weingarten canal
surfaces in IR$^{3}$ by using the definition of general canal surfaces.
During the study, we faced a very large equations. It was not possible to
give them all of course. So we had to make our processes via a computer time
to time.

Let $f$ and $g$ be smooth functions on a surface M in Euclidean 3-space E$%
^{3}$.The Jacobi function $\Phi \left( f,g\right) $ formed with $f,g$ is
defined by
\begin{equation*}
\Phi \left( f,g\right) =f_{s}g_{t}-f_{t}g_{s}
\end{equation*}%
where $f_{s}=\frac{\partial f}{\partial s}$ and $f_{t}=\frac{\partial f}{%
\partial t}$. In particular,a surface satisfying the Jacobi equation $\Phi ($%
K,H$)=0$ with respect to the Gaussian curvature K and the mean curvature H
on a surface M is called a Weingarten surface. Also, if a surface satisfies
a linear equation with respect to K and H, that is, aK+bH$=$c $\left( \text{%
a,b,c}\in \text{IR,(a,b,c)}\neq (\text{0,0,0})\right) $, then it is said to
be a linear Weingarten surface\cite{15}.

When the constant b=0, a linear Weingarten surface M reduces to a surface
with constant Gaussian curvature. When the constant a=0, a linear Weingarten
surface M reduces to a surface with constant mean curvature. In such a
sense, the linear Weingarten surfaces can be regarded as a natural
generalization of surfaces with constant Gaussian curvature or with constant
mean curvature\cite{15}.

If the second fundamental form of a surface M in E$^{3}$ is non-degenerate,
then it is regarded as a new pseudo-Riemannian metric. Therefore, the
Gaussian curvature K$_{\text{II}}$ of non-degenerate second fundamental form
can be defined formally on the Riemannian or pseudo-Riemannian manifold $($%
M,II$)$.We call the curvature K$_{\text{II}}$ the second Gaussian curvature
on M \cite{15}.

Following the Jacobi equation and the linear equation with respect to the
Gaussian curvature K, the mean curvature H and the second Gaussian curvature
K$_{\text{II}}$ an interesting geometric question is raised. Classify all
surfaces in Euclidean 3-space satisfying the conditions%
\begin{equation*}
\Phi (X,Y)=0
\end{equation*}%
\begin{equation*}
aX+bY=c
\end{equation*}%
where $\left( \text{X,Y}\right) \in \left\{ \text{K,H,K}_{\text{II}}\right\}
$ , X$\neq $Y and $(a,b,c)\neq (0,0,0).$ Let M be a surface immersed in
Euclidean 3-space, the first fundamental form of the surface M is defined by
\begin{equation*}
\text{I}=Edu^{2}+2Fdudv+Gdv^{2}
\end{equation*}%
where $E=<M_{s},M_{s}>$, $F=<M_{s},M_{t}>$, $G=<M_{t},M_{t}>$ are the
coefficients of I. A surface is called degenerate if it has the degenerate
first fundamental form. The second fundamental form of M is given by%
\begin{equation*}
\text{II}=edu^{2}+2fdudv+gdv^{2}
\end{equation*}%
where $e=\langle M_{ss},n\rangle $, $f=$ $\langle M_{st},n\rangle $, $%
g=\langle M_{tt},n\rangle $ and $n$ is the unit normal of M. The Gaussian
curvature K and the mean curvature H are given by, respectively%
\begin{equation}
\text{K}=\frac{eg-f^{2}}{EG-F^{2}},  \label{01}
\end{equation}%
\begin{equation}
\text{H}=\frac{Eg-2Ff+Ge}{2(EG-F^{2})}.  \label{02}
\end{equation}%
A regular surface is flat if and only if its Gaussian curvature vanishes
identically. A minimal surface in IR$^{3}$ is a regular surface for which
mean curvature vanishes identically\cite{Ga}.

Furthermore, the second Gaussian curvature K$_{\text{II}}$ of a surface is
defined by%
\begin{eqnarray}
\text{K}_{\text{II}} &=&\frac{1}{\left( \left\vert eg\right\vert
-f^{2}\right) ^{2}}\left\{ p-q\right\} .  \notag \\
&&  \label{03}
\end{eqnarray}%
where
\begin{equation*}
p=\left\vert
\begin{array}{ccc}
-\frac{1}{2}e_{tt}+f_{st}-\frac{1}{2}g_{ss} & \frac{1}{2}e_{s} & f_{s}-\frac{%
1}{2}e_{t} \\
f_{t}-\frac{1}{2}g_{s} & e & f \\
\frac{1}{2}g_{t} & f & g%
\end{array}%
\right\vert
\end{equation*}%
and
\begin{equation*}
q=\left\vert
\begin{array}{ccc}
0 & \frac{1}{2}e_{t} & \frac{1}{2}g_{s} \\
\frac{1}{2}e_{t} & e & f \\
\frac{1}{2}g_{s} & f & g%
\end{array}%
\right\vert
\end{equation*}%
A surface is called II-flat if the second Gaussian curvature vanishes
identically\cite{Yoon}. Having in mind the usual technique for computing the
second mean curvature by using the normal variation of the area functional
one gets%
\begin{equation}
\text{H}_{\text{II}}=\text{H}-\frac{1}{2\sqrt{\left\vert \det \text{II}%
\right\vert }}\sum_{i,j}\frac{\partial }{\partial u^{i}}\left( \sqrt{%
\left\vert \det \text{II}\right\vert }L^{ij}\frac{\partial }{\partial u^{j}}%
\left( \ln \sqrt{\left\vert K\right\vert }\right) \right)   \label{04}
\end{equation}%
where $u^{i}$ and $u^{j}$ stand for "$s$" and "$t$", respectively, and $%
\left( L^{ij}\right) =\left( L_{ij}\right) ^{-1},$ where $L_{ij}$ are the
coefficients of second fundamental forms\cite{Sv}. A surface is called
II-minimal if the second mean curvature vanishes identically\cite{Yoon}.

A canal surface is an envelope of a 1-parameter family of surface. The
envelope of a 1-parameter family $s\longrightarrow S^{2}\left( s\right) $ of
spheres in IR$^{3}$ is called a \textit{canal surface}\cite{Ga}. The curve
formed by the centers of the spheres \ is called \textit{center curve }of
the canal surface. The radius of canal surface is the function $r$ such that
$r(s)$ is the radius of the sphere $S^{2}\left( s\right) .$ Suppose that the
center curve of a canal surface is a unit speed curve $\alpha :I\rightarrow $%
IR$^{3}$. Then the canal surface can be parametrized by the formula
\begin{equation}
C\left( s,t\right) =\alpha \left( s\right) -R\left( s\right) T-Q\left(
s\right) \cos \left( t\right) N+Q\left( s\right) \sin \left( t\right) B
\label{0}
\end{equation}%
where%
\begin{gather}
R\left( s\right) =r(s)r^{\prime }(s)  \label{1} \\
Q\left( s\right) =\pm r(s)\sqrt{1-r^{\prime }(s)^{2}}.  \label{2}
\end{gather}%
All the tubes and the surfaces of revolution are subclass of the canal
surface.

\begin{thm}
The center curve of a canal surface M is a straight line if and only if M is
a surface of revolution for which no normal line to the surface is parallel
o the axis of revolution \cite{Ga}.
\end{thm}

\begin{thm}
The following conditions are equivalent for a canal surface M:

i. M is a tube parametrized by (\ref{0});

ii. the radius of M is constant;

iii. the radius vector of each sphere in family that defines the canal
surface M meets the center curve orthogonally \cite{Ga}.
\end{thm}

Coefficients of first and second fundamental forms of canal surface are%
\begin{eqnarray}
E\left( s,t\right) &=&Q^{2}\kappa ^{2}\cos ^{2}\left( t\right) +p_{1}\kappa
\cos \left( t\right) +2QR\kappa \tau \sin (t)+p_{2}  \notag \\
F\left( s,t\right) &=&-Q\left( R\kappa \sin (t)+Q\tau \right)  \label{3a} \\
G\left( s,t\right) &=&Q^{2}  \notag
\end{eqnarray}%
and%
\begin{eqnarray}
e\left( s,t\right) &=&\frac{-1}{r\left( s\right) }\left\{ E-Q\kappa \cos
\left( t\right) -p_{5}\right\}  \notag \\
f\left( s,t\right) &=&\frac{-1}{r\left( s\right) }F\left( s,t\right)
\label{3b} \\
g\left( s,t\right) &=&\frac{-1}{r\left( s\right) }G\left( s,t\right)  \notag
\end{eqnarray}%
Let us take $\psi \left( s,t\right) =\det I$ and $\phi \left( s,t\right)
=\det II.$ Thus, we have
\begin{equation}
\phi \left( s,t\right) =\frac{1}{r^{2}}\left\{ \psi \left( s,t\right)
-Q^{3}\kappa \cos \left( t\right) -Q^{2}p_{5}\right\}  \label{3c}
\end{equation}%
\begin{equation}
\psi \left( s,t\right) =Q^{2}\left\{
\begin{array}{c}
\kappa ^{2}\left( R^{2}+Q^{2}\right) \cos ^{2}\left( t\right) +\kappa
p_{1}\cos \left( t\right) \\
+1-2R^{\prime }+\left( R^{\prime }\right) ^{2}+\left( Q^{\prime }\right) ^{2}%
\end{array}%
\right\} .  \label{k2}
\end{equation}%
and%
\begin{eqnarray}
p_{1} &=&2\left( Q-QR^{\prime }+Q^{\prime }R\right)  \notag \\
p_{2} &=&Q^{2}\tau ^{2}+R^{2}\kappa ^{2}+\left( R^{\prime }\right)
^{2}+\left( Q^{\prime }\right) ^{2}-2R^{\prime }+1  \notag \\
p_{3} &=&p_{1}-Q  \label{3d} \\
p_{4} &=&p_{2}-p_{5}  \notag \\
p_{5} &=&\left( R^{\prime }\right) ^{2}+\left( Q^{\prime }\right)
^{2}-2R^{\prime }+1+RR^{\prime \prime }+QQ^{\prime \prime }  \notag
\end{eqnarray}%
If $Q(s)=0$, then the first and the second fundamental forms are degenerate.
So the canal surface is degenerate surface and the radius is $r(s)=\pm s+c$.
Furthermore, in the case $\kappa (s)=0$ and $1-2R^{\prime }+\left( R^{\prime
}\right) ^{2}+\left( Q^{\prime }\right) ^{2}=0,$ the radius is
\begin{equation*}
r(s)=\sqrt{s^{2}-2c_{1}s+2c_{2}}.
\end{equation*}%
Let the center curve be $\alpha \left( s\right) =\left( s,0,0\right) $. Then
$T=e_{1}$, $N=e_{2}$ and $B=e_{3}$. Hence, $R\left( s\right) =s-c_{1}$ and $%
C\left( s,t\right) $ is the curve in the plane $x=c_{1}$. The conditions
that $r(s)\neq \pm s+c$ and $\left( \kappa (s)=0,r(s)\neq \sqrt{%
s^{2}-2c_{1}s+2c_{2}}\right) $ are the necessary conditions to define a
non-degenerate canal surface with the equation (\ref{1}). At this point, we
can write the following theorem.

\begin{thm}
Let $M$ be a canal surface with the center curve $\alpha (s)$ and the radius
$r(s)$. If the center curve is a line then $M$ is a regular surface in IR$%
^{3}$ iff the radius is $r(s)\neq \pm s+c$ and $r(s)\neq \sqrt{%
s^{2}-2c_{1}s+2c_{2}}$.
\end{thm}

Additionally, if $\phi \left( s,t\right) =0$ then $M$ has degenerate second
fundamental form. A canal surface has degenerate second fundamental form if
canal surface is a surface of revolution with the radious $r(s)=\sqrt{%
s^{2}-2c_{1}s+2c_{2}}$ or $r(s)=c_{1}s+c_{2}$. From (\ref{01}), (\ref{02})
and (\ref{03}), we obtained the Gauss curvature, mean curvature such that%
\begin{eqnarray}
\text{K}(s,t) &=&\frac{-1}{\psi (s,t)r^{2}}\left\{ Q^{3}\kappa \cos \left(
t\right) +Q^{2}p_{5}-\psi (s,t)\right\}  \notag \\
&&  \label{k1}
\end{eqnarray}%
\begin{eqnarray}
\text{H}(s,t) &=&\frac{1}{2\psi (s,t)r^{2}}\left\{ Q^{3}\kappa \cos \left(
t\right) +Q^{2}p_{5}-2\psi (s,t)\right\} .  \notag \\
&&  \label{h1}
\end{eqnarray}

\section{Weingarten Type Canal Surfaces}

Let M be a canal surface with the center curve $\alpha (s)$ and the radius $%
r(s).$ The existence of a Weingarten relation $\Phi $(H,K) = 0 means that
curvatures H and K are functionally related, and since H and K are
differentiable functions depending on $s$ and $t$, this implies the Jacobian
condition $\Phi $(H,K)=0. More precisely the following condition%
\begin{equation}
\text{H}_{t}\left( \text{K}\right) _{s}\text{-H}_{s}\left( \text{K}\right)
_{t}=0  \label{w1}
\end{equation}%
needs to be satisfied. By using equations (\ref{k1}) and (\ref{h1}) we get%
\begin{equation}
\text{H}_{t}\left( \text{K}\right) _{s}\text{-H}_{s}\left( \text{K}\right)
_{t}=\frac{1}{2\psi ^{3}r^{4}}\sum\limits_{i=0}^{2}h_{i}\cos ^{i}\left(
t\right)   \label{w2}
\end{equation}%
where
\begin{eqnarray*}
h_{2} &=&-Q^{6}\kappa ^{2}\psi _{t}r^{\prime } \\
h_{1} &=&-Q^{2}\left\{
\begin{array}{c}
Q^{4}\kappa ^{2}\psi r^{\prime }\sin (t)+3Q\kappa \psi ^{2}\psi
_{t}r^{\prime }-Q\kappa ^{\prime }\psi ^{2}\psi _{t}r-3Q^{\prime }\kappa
\psi ^{2}\psi _{t}r \\
+Q\kappa \psi \theta _{t}r^{\prime }-2Q\kappa \psi _{t}\theta r^{\prime }%
\end{array}%
\right\}  \\
h_{0} &=&\psi \theta _{t}\left\{ r\psi \psi _{s}-r^{\prime }(\psi
^{2}-\theta )\right\} +\psi _{t}\left\{ 3\theta r^{\prime }(\psi ^{2}-\frac{1%
}{3}\theta )-r\psi ^{2}\theta _{s}\right\}  \\
&&+Q^{3}\psi \kappa \left\{ r\psi \psi _{s}-r^{\prime }(\psi ^{2}-\theta
)\right\} \sin (t).
\end{eqnarray*}%
and $\theta =\psi -Q^{2}p_{5}.$ The Jacobian condition requires $%
h_{0}=h_{1}=h_{2}=0.$ From $h_{2}=0,$ the cases $\kappa =0,\psi
_{t}=0,r^{\prime }=0$ are possible$.$ If $\kappa =0$ then, $\psi =Q^{2}p_{2}$
and $h_{0}=h_{1}=0$ satisfies. If $\psi _{t}=0$ then, from (\ref{k2}) $%
\kappa =0$. If $r=c$ then, from $h_{1}=0$ we have
\begin{equation*}
c^{3}\psi _{t}\kappa ^{\prime }=0
\end{equation*}%
In this case, if $\psi _{t}=0$ then, from (\ref{k2}) $\kappa =0$ and so (\ref%
{k1}) satisfy. If $\kappa =a=$constant then, $h_{1}=h_{2}=0$ and $h_{0}=0$
is
\begin{equation*}
\psi _{s}(ac^{3}\sin (t)+\theta _{t})-\psi _{t}\theta _{s}=0.
\end{equation*}%
Since the relation $ac^{3}\sin (t)+\theta _{t}\neq 0$ and $\psi _{t}\neq 0$
then, $\psi _{s}=0$ and $\psi _{t}\theta _{s}=0.$ Also from (\ref{3d}), $%
p_{1}=2c,$ $p_{5}=1$ and $p_{2}=c^{2}\tau ^{2}+1.$ Thus $\psi _{s}=\theta
_{s}=0$ satisfy. Hence we proved the following theorem.

\begin{thm}
\bigskip Let M be a regular canal surface then, M is a (H,K)-Weingarten
canal surface if M is one of the surfaces that surface of revolution,
cyclinder and tubular surface whose centered curve with nonzero constant
curvature.
\end{thm}

Thus, there are three cases for (\ref{0}) such as $\left( r^{\prime }\neq
0,\kappa =0\right) ,$ $\left( r^{\prime }=0,\kappa \neq 0\right) $ and $%
\left( r^{\prime }=0,\kappa =0\right) $, If $\kappa =0$ then, all of the
coefficients $h_{i}$ in (\ref{w2}) are zero so the condition (\ref{w1}) is
satisfy. If $r^{\prime }=0$ and $\kappa \neq 0$ then $r(s)=c\neq 0$ and (\ref%
{0}) turns to a tubular suface such that%
\begin{equation}
C\left( s,t\right) =\alpha \left( s\right) \mp c\cos \left( t\right) N\pm
c\sin \left( t\right) B.  \label{w4}
\end{equation}%
If $\kappa =0$ then let assume that the center curve is the $x-$axis (\ref{0}%
) turns to a surface of revolution and a cylinder such that%
\begin{equation}
C\left( s,t\right) =\left( s-r(s)r^{\prime }(s),\mp r(s)\sqrt{1-r^{\prime
}(s)^{2}}\cos \left( t\right) ,\pm r(s)\sqrt{1-r^{\prime }(s)^{2}}\sin
\left( t\right) \right)   \label{w6}
\end{equation}%
and%
\begin{equation}
C\left( s,t\right) =\left( s,\mp c\cos \left( t\right) ,\pm c\sin \left(
t\right) \right)   \label{w5}
\end{equation}%
respectively.

\bigskip From (\ref{03}), we can write the term $p$ as
\begin{equation}
p=\left( -\frac{1}{2}e_{tt}+f_{st}-\frac{1}{2}g_{ss}\right) \phi +\left(
f_{t}-\frac{1}{2}g_{s}\right) \left\{ \left( f_{s}-\frac{1}{2}e_{t}\right) f-%
\frac{1}{2}e_{s}g\right\}   \label{wa}
\end{equation}%
by taking $g_{t}=0.$ From  (\ref{h1}), (\ref{wa}) and with the aid of prog
2, the Jacobi function $\Phi $(H,K$_{\text{II}}$) is obtained a polynomial
expressions in cos$\left( t\right) $ such that%
\begin{equation}
\text{H}_{t}\left( \text{K}_{\text{II}}\right) _{s}-\text{H}_{s}\left( \text{%
K}_{\text{II}}\right) _{t}=\frac{1}{denom}\sum\limits_{i=0}^{6}g_{i}\cos
^{i}\left( t\right) .  \label{w7}
\end{equation}%
For n=6 in Prog.2, $g_{6}$ is
\begin{equation*}
g_{6}=-Q^{13}\kappa ^{5}r^{4}\psi _{t}\left\{ 4Q\kappa r^{\prime }-3Q\kappa
^{\prime }r-5Q^{\prime }\kappa r\right\} .
\end{equation*}%
The Jacobian condition $\Phi $(H,K$_{\text{II}}$)=0 requires $%
g_{0}=g_{1}=...g_{6}=0.$ In the case $\kappa =0,$ $g_{1}=...g_{6}=0$
satisfies. By using prog.1 for $\kappa =\tau =0,$ it is easy to see that $%
g_{0}=0$ satisfy. If $\psi _{t}=0$ then, from (\ref{3a}) we obtain $\kappa =0
$ so $\Phi $(H,K$_{\text{II}}$)=0 satisfy also. If $4Q\kappa r^{\prime
}-3Q\kappa ^{\prime }r-5Q^{\prime }\kappa r=0$ then, from (\ref{1}) and (\ref%
{2}), we obtain
\begin{equation}
5\kappa rr^{\prime }r^{\prime \prime }+(\left( r^{\prime }\right)
^{2}-1)(3r\kappa ^{\prime }+r^{\prime }\kappa )=0.  \label{w77}
\end{equation}%
It may be hard to solve (\ref{w77}), but ofcourse we consider the special
solutions of (\ref{w77}). If $r$ is constant then from (\ref{w77}) $\kappa $
is non zero constant. If $r=c_{1}s+c_{2\text{ }}$then, (\ref{w77}) turn to%
\begin{equation*}
(\left( c_{1}\right) ^{2}-1)(3r\kappa ^{\prime }+c_{1}\kappa )=0
\end{equation*}%
and the solutions are $c_{1}=\pm 1$(M is not regular) or
\begin{equation*}
\kappa =\frac{c_{3}}{r^{1/3}}.
\end{equation*}%
If $\kappa =c_{1}$ is non-zero constant in (\ref{w77}) then, (\ref{w77})
turn to%
\begin{equation*}
rr^{\prime }(5rr^{\prime \prime }+(r^{\prime })^{2}-1)=0
\end{equation*}%
and the real solution is $r=$constant$\neq \pm 1.$ Thus we can give the
following theorem.

\begin{thm}
Let M be a regular (H,K$_{II}$)-Weingarten canal surface then followings are
ture for M.

i.M is the surface of revolution,

ii.M is a canal surface with $r=c_{1}s+c_{2\text{ }}$,($c_{1}\neq \pm $1)
and with the centered curve whose curvature is $\kappa =\frac{c_{3}}{r^{1/3}}
$,

iii.M is a tubular surface whose centered curve with non-zeroconstant
curvature..
\end{thm}

Jacobi function $\Phi $(K,K$_{\text{II}}$) is obtained a polynomial
expressions in cos$\left( t\right) $ by using (\ref{k1}) and (\ref{wa}) as
follows.%
\begin{equation}
\text{K}_{t}\left( \text{K}_{\text{II}}\right) _{s}-\text{K}_{s}\left( \text{%
K}_{\text{II}}\right) _{t}=\frac{1}{denom}\sum\limits_{i=0}^{6}f_{i}\cos
^{i}\left( t\right)   \label{w8}
\end{equation}%
Jacobian condition $\Phi $(K,K$_{\text{II}}$)=0 requires $%
f_{0}=f_{1}=...f_{6}=0.$For $i=6$ in Prog.3, $f_{6}$ is the same as $g_{6}$
in (\ref{w7}).
\begin{equation*}
f_{6}=Q^{13}\kappa ^{5}r^{4}\psi _{t}\left\{ 4Q\kappa r^{\prime }-3Q\kappa
^{\prime }r-5Q^{\prime }\kappa r\right\} .
\end{equation*}%
If $\kappa =0$ then, $\Phi $(K,K$_{\text{II}}$)=0 satisfy. If $\psi _{t}=0$
then, from (\ref{k2}) $\kappa =0$ and so $\Phi $(K,K$_{\text{II}}$)=0
satisfy. If $4Q\kappa r^{\prime }-3Q\kappa ^{\prime }r-5Q^{\prime }\kappa r=0
$ then, we obtain the same differential equation (\ref{w77}) and we consider
the special solution again. If $r$ is constant then from (\ref{w77}) $\kappa
$ is non-zero constant. In this case, using Prog.1 and Prog.3, we obtain $%
f_{6}=f_{5}=f_{4}=0$ and $f_{3}=-3\kappa ^{5}r^{2}\tau ^{\prime }\sin ^{2}(t)
$, $f_{2}=-6\kappa ^{4}r\tau ^{\prime }\sin ^{2}(t),$ $f_{1}=\kappa
^{3}(4\kappa ^{2}r^{2}-1)\tau ^{\prime }\sin ^{2}(t),$ $f_{0}=2\kappa
^{4}r\tau ^{\prime }\sin ^{2}(t).$ The Jacobien condition reguires that $%
\tau $ is constant. If $r=c_{1}s+c_{2\text{ }}$then, $\kappa =\frac{c_{3}}{%
r^{1/3}}$ and also, $\Phi $(K,K$_{\text{II}}$)=0 satisfy. Thus, we can write
the following theorem.

\begin{thm}
Let M be a regular (K,K$_{II}$)-Weingarten canal surface then followings are
ture for M.

i.M is the surface of revolution,

ii.M is a canal surface with $r=c_{1}s+c_{2\text{ }}$,($c_{1}\neq \pm $1)
and with the centered curve whose curvature is $\kappa =\frac{c_{3}}{r^{1/3}}
$,

iii.M is a tubular surface whose centered curve is cylindirical helix.
\end{thm}

\section{Linear\textbf{\ }Weingarten type Canal Surfaces}

Let M be a canal surface with the center curve $\alpha (s)$ and the radius $%
r(s)$ then M is called $\left( \text{K,H}\right) -$linear Weingarten surface
if Gaussian and the mean curvatures of M satisfies a linear equation with
the constants a, b and d such that \
\begin{equation*}
\text{aK+bH}=\text{d}.
\end{equation*}%
By using equations (\ref{k1}) and (\ref{h1}), we get the relation between K$%
\left( s,t\right) $ and H$\left( s,t\right) $ such as
\begin{equation}
\text{H}\left( s,t\right) \text{+}\frac{1}{2}\text{K}\left( s,t\right) =-%
\frac{1}{2r^{2}}.  \label{30}
\end{equation}%
Thus we have the following theorem.

\begin{thm}
Let M be a regular canal surface then M is $\left( \text{K,H}\right) -$%
linear Weingarten surface if and only if M is a tubular surface.
\end{thm}

From\ (\ref{04}),we can write
\begin{equation}
denom(\delta )\sqrt{\mu _{1}\phi }\text{H}_{\text{II}}-denom(\delta )\sqrt{%
\mu _{1}\phi }\text{H}=\frac{\text{numer(}\delta )}{2}  \label{3}
\end{equation}%
where $\delta =\delta _{1}+\delta _{2},$%
\begin{eqnarray*}
\delta _{1} &=&\frac{\partial }{\partial s}\left( \sqrt{\mu _{1}\phi }L^{11}%
\frac{\partial }{\partial s}\left( \ln \sqrt{\mu _{2}K}\right) +\sqrt{\mu
_{1}\phi }L^{12}\frac{\partial }{\partial t}\left( \ln \sqrt{\mu _{2}K}%
\right) \right)  \\
\delta _{2} &=&\frac{\partial }{\partial t}\left( \sqrt{\mu _{1}\phi }L^{21}%
\frac{\partial }{\partial s}\left( \ln \sqrt{\mu _{2}K}\right) +\sqrt{\mu
_{1}\phi }L^{22}\frac{\partial }{\partial t}\left( \ln \sqrt{\mu _{2}K}%
\right) \right) ,
\end{eqnarray*}%
\begin{equation*}
\mu _{1}=\left\{
\begin{array}{c}
1 \\
-1%
\end{array}%
\begin{array}{c}
;\text{if }\phi >0 \\
;\text{if }\phi <0%
\end{array}%
\right. \text{ \ and }\mu _{2}=\left\{
\begin{array}{c}
1 \\
-1%
\end{array}%
\begin{array}{c}
;\text{if }K>0 \\
;\text{if }K<0%
\end{array}%
\right.
\end{equation*}%
and with the aid of Prog.4
\begin{equation*}
denom(\delta )=4\phi ^{2}\psi ^{2}r^{2}\sqrt{\mu _{1}\phi }.
\end{equation*}%
Thus, (\ref{3}) turn to%
\begin{equation}
4\phi ^{3}\psi ^{2}r^{2}\text{H}_{\text{II}}-4\phi ^{3}\psi ^{2}r^{2}\text{H}%
=\frac{\text{numer(}\delta )}{2}  \label{3aa1}
\end{equation}%
In this case, from (\ref{3aa1}), if $\phi ^{3}\psi ^{2}r^{2}$ and numer($%
\delta )$ are constant then, we can say that there is a linear relation
between H$_{II}$ and H. By using (\ref{3c}), (\ref{k2}) and Prog.5
\begin{equation}
(r^{2}\phi ^{3}\psi ^{2})_{t}=\frac{1}{r^{4}}\sin
(x)\sum\limits_{i=0}^{9}m_{i}\cos ^{i}(x)  \label{3aa}
\end{equation}%
for $i=9$
\begin{equation*}
m_{9}=-10Q^{10}\kappa ^{10}(Q^{2}+R^{2})^{5}
\end{equation*}%
from $m_{9}=0$ then $\kappa =0$ and also all of $m_{i}$ are zero. In the
case of $\kappa =0,$ $r^{2}\phi ^{3}\psi ^{2}$ is%
\begin{equation}
r^{2}\phi ^{3}\psi ^{2}=\frac{1}{r^{4}}Q^{10}(p_{2}-p_{5})^{3}(p_{2})^{2}.
\label{3a0}
\end{equation}%
By using prog.5 and 6, the real non-zero solutions of $(r^{2}\phi ^{3}\psi
^{2})_{s}=0$ are $r=\pm \sqrt{s^{2}-2c_{1}s+2c_{2}}$ and $r=c_{1}s+c_{2}.$
For first $r,$ M is degenerate, and for the second $r$, $\phi =0$. \ Thus we
have the following theorem.

\begin{thm}
Let M be a regular canal surface then, there is \ no (H,H$_{\text{II}}$%
)-linear Weingarten surface in IR$^{3}.$
\end{thm}

From\ (\ref{01}),
\begin{equation*}
\phi \left( s,t\right) =\psi \left( s,t\right) K\left( s,t\right) ,
\end{equation*}%
and by using (\ref{3c}), (\ref{k1}) and (\ref{h1}), we get
\begin{equation}
\phi ^{2}\text{K}_{\text{II}}-A\psi \text{K}=(B-q)  \label{3a1}
\end{equation}%
where
\begin{equation*}
A=\left( -\frac{1}{2}e_{tt}+f_{st}-\frac{1}{2}g_{ss}\right)
\end{equation*}%
and
\begin{equation*}
B=\left( f_{t}-\frac{1}{2}g_{s}\right) \left\{ \left( f_{s}-\frac{1}{2}%
e_{t}\right) f-\frac{1}{2}e_{s}g\right\} .
\end{equation*}%
If $\phi ,$ $r^{4}A\psi $ and $r^{4}(B-q)$ are nonzero constans then, we
called M is (K$_{\text{II}},$K)-linear Weingarten canal surface.%
\begin{equation*}
\frac{\partial \phi }{\partial t}=0
\end{equation*}%
\begin{equation*}
\frac{\kappa Q^{2}}{r^{4}}\left\{ -2\kappa \left( R^{2}+Q^{2}\right) \cos
\left( t\right) -\left( p_{1}-Q\right) \right\} \sin \left( t\right) =0
\end{equation*}%
Thus, $\kappa =0$ and for $\kappa =0,$%
\begin{equation*}
\phi =\frac{-Q^{2}}{r^{4}}\left\{ RR^{\prime \prime }+QQ^{\prime \prime
}\right\}
\end{equation*}%
and by using (\ref{01}) and (\ref{02}), we obtain%
\begin{equation*}
\frac{-Q^{2}}{r^{4}}\left\{ RR^{\prime \prime }+QQ^{\prime \prime }\right\}
=rr^{\prime \prime }((r^{\prime })^{2}+rr^{\prime \prime }-1)
\end{equation*}%
since $\frac{\partial \phi }{\partial s}=0$ then, $\frac{-Q^{2}}{r^{4}}%
\left\{ RR^{\prime \prime }+QQ^{\prime \prime }\right\} =c=$constant so we
can write
\begin{equation*}
rr^{\prime \prime }((r^{\prime })^{2}+rr^{\prime \prime }-1)=c
\end{equation*}%
There are the only three real non-zero solution of last equation for $c=0$
such that, $r=\pm \sqrt{s^{2}-2c_{1}s+2c_{2}}$ and $r=c_{1}s+c_{2}$ but in
the case of $r=\pm \sqrt{s^{2}-2c_{1}s+2c_{2}}$ and $r=c_{1}s+c_{2}$, M have
degenerate first and second fundamental forms. Thus, we can give the
following theorem.

\begin{thm}
Let M be a regular canal surface then, there is \ no (K,K$_{\text{II}}$%
)-linear Weingarten surface in IR$^{3}.$
\end{thm}

By substituting H in the equation (\ref{30}) into (\ref{3a1}), we get%
\begin{equation}
8\phi ^{3}\psi ^{2}r^{2}\text{H}_{\text{II}}+4\phi ^{3}\psi ^{2}r^{2}\text{K}%
=\text{numer(}\delta )-\frac{4\phi ^{3}\psi ^{2}r^{2}}{r^{2}}  \label{3a2}
\end{equation}%
and also, we found before that $\phi =0$ when $\phi ^{3}\psi ^{2}r^{2}$ is a
constant then, we have the following theorem.

\begin{thm}
Let M be a regular canal surface then, there is \ no (K,H$_{\text{II}}$%
)-linear Weingarten surface in IR$^{3}.$
\end{thm}

Similarly, by substituting K in the equation (\ref{30}) into (\ref{3a1}), we
get
\begin{equation}
\phi ^{2}r^{2}\text{K}_{\text{II}}+2r^{2}A\psi \text{H}=r^{2}(B-q)-A\psi .
\label{3aa0}
\end{equation}%
Also, we found before that M have degenerate first and second fundamental
forms when $\left\{ \psi \left( s,t\right) -Q^{3}\kappa \cos \left( t\right)
-Q^{2}p_{5}\right\} $ is a constant then, we have the following theorem.

\begin{thm}
Let M be a regular canal surface then, there is \ no (H,K$_{\text{II}}$%
)-linear Weingarten surface in IR$^{3}.$
\end{thm}

We can find the relation between K$_{\text{II}}$ and H$_{\text{II}}$ by
using (\ref{3a1}) and (\ref{3a2}) as follow
\begin{equation}
8\phi ^{3}\psi ^{2}r^{2}\text{H}_{\text{II}}+\frac{4\phi ^{3}\psi \left\{
\phi r^{2}\right\} ^{2}}{r^{2}A}\text{K}_{\text{II}}=\text{numer(}\delta
)+\Gamma  \label{3a3}
\end{equation}%
where $\Gamma =\frac{4\phi ^{3}\psi (r^{2}(B-q)-A\psi )}{A}.$ Since $\phi =0$
when $\phi ^{3}\psi ^{2}r^{2}$ is a constant then, we give the following
theorem.

\begin{thm}
Let M be a regular canal surface then, there is \ no (K$_{\text{II}}$,H$_{%
\text{II}}$)-linear Weingarten surface in IR$^{3}.$
\end{thm}

It easy to obtain the relations of trible of $\left\{ \text{H},\text{K},%
\text{H}_{\text{II}},\text{K}_{\text{II}}\right\} .$ From (\ref{30}) and (%
\ref{3aa1}),
\begin{equation}
2(1-4\gamma _{1})\text{H}+\text{K}+8\gamma _{1}\text{H}_{\text{II}}=\frac{%
r^{2}\text{numer(}\delta )-1}{r^{2}},  \label{3a4}
\end{equation}%
from (\ref{30}) and (\ref{3aa1}), \
\begin{equation}
2\text{H}+(1-2r^{4}A\psi )\text{K}+2(\gamma _{2})^{2}\text{K}_{\text{II}}=%
\frac{2r^{6}(B-q)-1}{r^{2}},  \label{3a5}
\end{equation}%
from (\ref{3a1}) and (\ref{3a2}),%
\begin{equation}
\psi r^{2}(4\phi ^{3}\psi -r^{2}A)\text{K}+8\gamma _{1}\text{H}_{\text{II}%
}+(\gamma _{2})^{2}\text{K}_{\text{II}}=\text{numer}(\delta
)+r^{4}(B-q)-4\phi ^{3}\psi ^{2},  \label{3a6}
\end{equation}%
from (\ref{3aa0}) and (\ref{3a3}),%
\begin{equation}
2r^{4}A\psi \text{H+}8\gamma _{1}\text{H}_{\text{II}}+\frac{(4\phi ^{3}\psi
+r^{2}A)(\gamma _{2})^{2}}{r^{2}A}\text{K}_{\text{II}}=\text{numer(}\delta
)+\gamma _{3}  \label{3a7}
\end{equation}%
where%
\begin{eqnarray*}
\gamma _{1} &=&\phi ^{3}\psi ^{2}r^{2} \\
\gamma _{2} &=&\phi r^{2} \\
\gamma _{3} &=&\frac{(4\phi ^{3}\psi +Ar^{2})(r^{2}(B-q)-A\psi )}{A}.
\end{eqnarray*}%
In (\ref{3a4}), (\ref{3a5}), (\ref{3a6}) and (\ref{3a7}), the conditions $%
\gamma _{1}$ and $\gamma _{2}$ are constants requires the first and the
second fundamental forms are degenerate. Thus, we have the following theorem.

\begin{thm}
Let M be a regular canal surface then, there is \ no (H,K,H$_{\text{II}}$),
(H,K,K$_{\text{II}}$), (K,K$_{\text{II}}$,H$_{\text{II}}$) and (H,K$_{\text{%
II}}$,H$_{\text{II}}$)-linear Weingarten surfaces in IR$^{3}.$
\end{thm}

Finally, from (\ref{3a4}) and (\ref{3a5}) we get%
\begin{equation*}
4(\text{1}-2\gamma _{1}))\text{H}+2(1-r^{4}A\psi )\text{K}+8\gamma _{1}\text{%
H}_{\text{II}}+2(\gamma _{2})^{2}\text{K}_{\text{II}}=\frac{r^{2}\text{numer(%
}\delta )+2r^{6}(B-q)-2}{r^{2}}
\end{equation*}%
and the condition $\gamma _{1}$ is constant requires the second fundamental
form is degenerate. Thus, we have the following theorem.

\begin{thm}
Let M be a regular canal surface then, there is \ no (H,K,H$_{\text{II}}$,K$%
_{\text{II}}$)-linear Weingarten surfaces in IR$^{3}.$
\end{thm}

\textbf{Prog.1}

R(s):=r(s)*diff(r(s),s):

Q(s):=r(s)*((1-diff(r(s),s)\symbol{94}2)\symbol{94}(1/2)):

G(s):=(Q(s))\symbol{94}2:

p1(s):=2*Q(s)+2*R(s)*diff(Q(s),s)-2*Q(s)*diff(R(s),s):

p2(s):=(Q(s)\symbol{94}2)*(tau(s)\symbol{94}2)+(R(s)\symbol{94}2)*(kappa(s)%
\symbol{94}2)

+diff(R(s),s)\symbol{94}2+diff(Q(s),s)\symbol{94}2-2*diff(R(s),s)+1;

p3(s):=p1(s)-Q(s):

p5(s):=diff(R(s),s)\symbol{94}2+diff(Q(s),s)\symbol{94}2-2*diff(R(s),s)+1

+R(s)*diff(diff(R(s),s),s)+Q(s)*diff(diff(Q(s),s),s);

p4(s):=p2(s)-p5(s):

E(s,t):=(Q(s)\symbol{94}2)*(kappa(s)\symbol{94}2)*(cos(t))\symbol{94}%
2+p1(s)*kappa(s)*cos(t)

+2*Q(s)*R(s)*kappa(s)*tau(s)*sin(t)+p2(s);

F(s,t):=-Q(s)*R(s)*kappa(s)*sin(t)-G(s)*tau(s):

psi(s,t):=E(s,t)*G(s)-(F(s,t))\symbol{94}2:

theta(s,t):=-(Q(s)\symbol{94}2)*p5(s)+psi(s,t):

\textbf{Prog.2}

G(s):=(Q(s))\symbol{94}2;

e:=(-1/r(s))*(E(s,t)-Q(s)*kappa(s)*cos(t)-p5(s)):

f:=(-1/r(s))*F(s,t):

g:=(-1/r(s))*G(s):

es:=diff(e,s):

fs:=diff(f,s):

gs:=diff(g,s):

ess:=diff(diff(e,s),s):

fss:=diff(diff(f,s),s):

gss:=diff(diff(g,s),s):

et:=diff(e,t):

ft:=diff(f,t):

gt:=diff(g,t):

ett:=diff(diff(e,t),t):

ftt:=diff(diff(f,t),t):

gtt:=diff(diff(g,t),t):

est:=diff(diff(e,s),t):

fst:=diff(diff(f,s),t):

gst:=diff(diff(g,s),t):

ets:=diff(diff(e,t),s):

fts:=diff(diff(f,t),s):

gts:=diff(diff(g,t),s):

phi(s,t):=(1/(r(s))\symbol{94}2)*(psi(s,t)-(Q(s)\symbol{94}%
3)*kappa(s)*cos(t)-(Q(s)\symbol{94}2)*p5(s)):

V1:=Matrix([[phi(s,t)*((-ett/2)+fst-(gss/2)),0],[0,phi(s,t)*((-ett/2)+fst

\ \ \ \ \ \ -(gss/2))]])+Matrix([[ft-(gs/2),((fs-(et/2))*e-(f*es/2))],[(gt/2)

\ \ \ \ \ ,(fs-(et/2))*f-(g*es/2)]]):

V2:=Matrix([[0,et/2,gs/2],[et/2,e,f],[gs/2,f,g]]):

v1:=LinearAlgebra:-Determinant(V1):

v2:=LinearAlgebra:-Determinant(V2):

K2:=simplify((v1-v2)/(phi(s,t))):

H:=-1/2*(-Q(s)\symbol{94}3*kappa(s)*cos(t)+theta(s,t)+psi(s,t)\symbol{94}%
2)/psi(s,t)/r(s):

simplify(coeff(numer(subs(sin(t)=B,subs(cos(t)=A,simplify(diff(H,t)

*diff(K2,s)-diff(H,s)*diff(K2,t))))),A,n),'size');

\textbf{Prog.3}

G(s):=(Q(s))\symbol{94}2;

e:=(-1/r(s))*(E(s,t)-Q(s)*kappa(s)*cos(t)-p5(s)):

f:=(-1/r(s))*F(s,t):

g:=(-1/r(s))*G(s):

es:=diff(e,s):

fs:=diff(f,s):

gs:=diff(g,s):

ess:=diff(diff(e,s),s):

fss:=diff(diff(f,s),s):

gss:=diff(diff(g,s),s):

et:=diff(e,t):

ft:=diff(f,t):

gt:=diff(g,t):

ett:=diff(diff(e,t),t):

ftt:=diff(diff(f,t),t):

gtt:=diff(diff(g,t),t):

est:=diff(diff(e,s),t):

fst:=diff(diff(f,s),t):

gst:=diff(diff(g,s),t):

ets:=diff(diff(e,t),s):

fts:=diff(diff(f,t),s):

gts:=diff(diff(g,t),s):

phi(s,t):=(1/(r(s))\symbol{94}2)*(psi(s,t)-(Q(s)\symbol{94}%
3)*kappa(s)*cos(t)-(Q(s)\symbol{94}2)*p5(s)):

V1 := Matrix([[phi(s,t)*((-ett/2)+fst-(gss/2)),0],[0,phi(s,t)*((-ett/2)+fst

\ \ \ \ \ -(gss/2))]])+
Matrix([[ft-(gs/2),((fs-(et/2))*e-(f*es/2))],[(gt/2),(fs

\ \ \ \ \ -(et/2))*f-(g*es/2)]]):

V2 := Matrix([[0,et/2,gs/2],[et/2,e,f],[gs/2,f,g]]):

v1 := LinearAlgebra:-Determinant(V1):

v2 := LinearAlgebra:-Determinant(V2):

K2:=simplify((v1-v2)/(phi(s,t))):

K:=1/psi(s,t)*(-Q(s)\symbol{94}3*kappa(s)*cos(t)+theta(s,t))/r(s)\symbol{94}%
2:

simplify(coeff(numer(subs(sin(t)=B,subs(cos(t)=A,simplify(diff(K,t)

*diff(K2,s)-diff(K,s)*diff(K2,t))))),A,i),'size');

\textbf{Prog.4}

G(s):=(Q(s))\symbol{94}2:

e:=(-1/r(s))*(E(s,t)-Q(s)*kappa(s)*cos(t)-p5(s)):

f:=(-1/r(s))*F(s,t):

g:=(-1/r(s))*G(s):

K:=phi(s,t)/psi(s,t);

L11:=g/(phi(s,t)):

L12:=-f/(phi(s,t)):

L21:=-f/(phi(s,t)):

L22:=e/(phi(s,t)):

delta1:=simplify(diff(((mu1*(phi(s,t)))\symbol{94}(1/2))*L11*diff(ln((mu2*K)%
\symbol{94}(1/2)),s)

\ \ \ \ \ \ \ \ \ +((mu1*(phi(s,t)))\symbol{94}(1/2))*L12*diff(ln((mu2*K)%
\symbol{94}(1/2)),t),s)):

delta2:=simplify(diff(((mu1*(phi(s,t)))\symbol{94}(1/2))*L21*diff(ln((mu2*K)%
\symbol{94}(1/2)),s)

\ \ \ \ \ \ \ \ \ +((mu1*(phi(s,t)))\symbol{94}(1/2))*L22*diff(ln((mu2*K)%
\symbol{94}(1/2)),t),t)):

denom(delta1+delta2);

\textbf{Prog.5}

E(s,t):=(Q(s)\symbol{94}2)*(kappa(s)\symbol{94}2)*(cos(t))\symbol{94}%
2+p1(s)*kappa(s)*cos(t)

\ \ \ \ \ \ \ \ \ +2*Q(s)*R(s)*kappa(s)*tau(s)*sin(t)+p2(s):

F(s,t):=-Q(s)*R(s)*kappa(s)*sin(t)-G(s)*tau(s):

G(s):=(Q(s))\symbol{94}2:

psi(s,t):=E(s,t)*G(s)-F(s,t)\symbol{94}2:

phi(s,t):=simplify((1/(r(s))\symbol{94}2)*(psi(s,t)-(Q(s)\symbol{94}%
3)*kappa(s)*cos(t)

\ \ \ \ \ \ \ \ \ \ -(Q(s)\symbol{94}2)*p5(s))):

simplify(coeff(numer(subs(sin(t)=B,subs(cos(t)=A,

diff(expand(simplify(phi(s,t)\symbol{94}3*psi(s,t)\symbol{94}2*r(s)\symbol{94%
}2)),t)))),A,9),'size');

\textbf{Prog.6}

dsolve(diff(expand(simplify(phi(s,t)\symbol{94}3*psi(s,t)\symbol{94}2*r(s)%
\symbol{94}2)),s), \{ r(s) \} );



\end{document}